\documentclass{amsart}
\usepackage[dvips]{epsfig}
\newtheorem{thm}{Theorem}[subsection]
 \newtheorem{cor}[thm]{Corollary}
 \newtheorem{lem}[thm]{Lemma}
 \newtheorem{prop}[thm]{Proposition}
 \theoremstyle{definition}

 \numberwithin{equation}{subsection}

\begin{document}
\title[Positive Dehn Twist Expressions for some New Involutions II]
 {Positive Dehn Twist Expressions for some New Involutions in the Mapping Class Group II\\}

\author{ Yusuf Z. Gurtas }

\address{Department of Mathematics, Suffolk CCC, Selden, NY, USA}

\email{gurtasy@sunysuffolk.edu}

\thanks{}




\dedicatory{}

\commby{}

\begin{abstract}
This article is continuation from \cite{Gu}. The positive Dehn
twist expressions for the generalization of the involutions
described in \cite{Gu} are presented. The homeomorphism types of
the Lefschetz fibrations they define are determined for several
examples.

\end{abstract}

\maketitle

\section*{Introduction}

In \cite{Gu} the author presented the positive Dehn twist
expression for a new set of involutions that are obtained by
combining two well known involutions in the mapping class group
$M_g$ of a $2$-dimensional, closed, compact, oriented surface
$\Sigma_{g}$ of genus $g>0$, one of which is the hyperelliptic
involution, Figure \ref{twoinvolutions.fig}. One can extend these
new involutions by gluing them together. It is the purpose of this
article to find the positive Dehn twist expressions for these
extended involutions and compute the signatures of the symplectic
Lefschetz fibrations that they describe.

\section{Review of the Simple Case}

\begin{figure}[htbp]
     \centering  \leavevmode
     \psfig{file=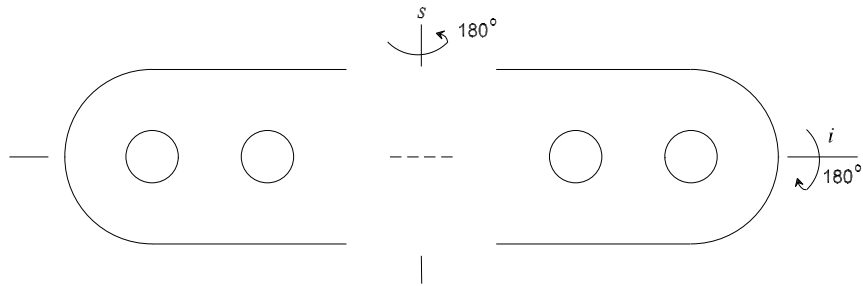,width=4.50in,clip=}
     \caption{The vertical and horizontal involutions}
     \label{twoinvolutions.fig}
 \end{figure}

Let $i$ represent the hyperelliptic (horizontal) involution and
$s$ represent the vertical involution as shown in Figure
\ref{twoinvolutions.fig}.

\begin{figure}[htbp]
     \centering  \leavevmode
     \psfig{file=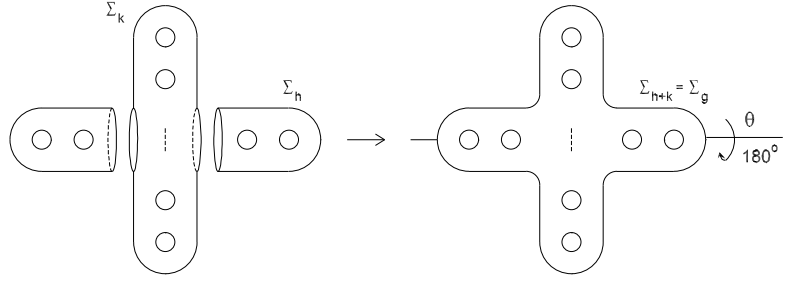,width=4.50in,clip=}
     \caption{The involution $\theta$ on the surface $\Sigma_{h+k}$}
     \label{gluingtwoinvolutions.fig}
 \end{figure}

If $i$ is the horizontal involution on a surface $\Sigma_{h}$ and
$s$ is the vertical involution on a surface $\Sigma_{k}$,
$k-$even, then let $\theta$ be the horizontal involution on the
surface $\Sigma_{g}$, where $g=h+k$, obtained as in Figure
\ref{gluingtwoinvolutions.fig}.

\begin{figure}[htbp]
     \centering  \leavevmode
     \psfig{file=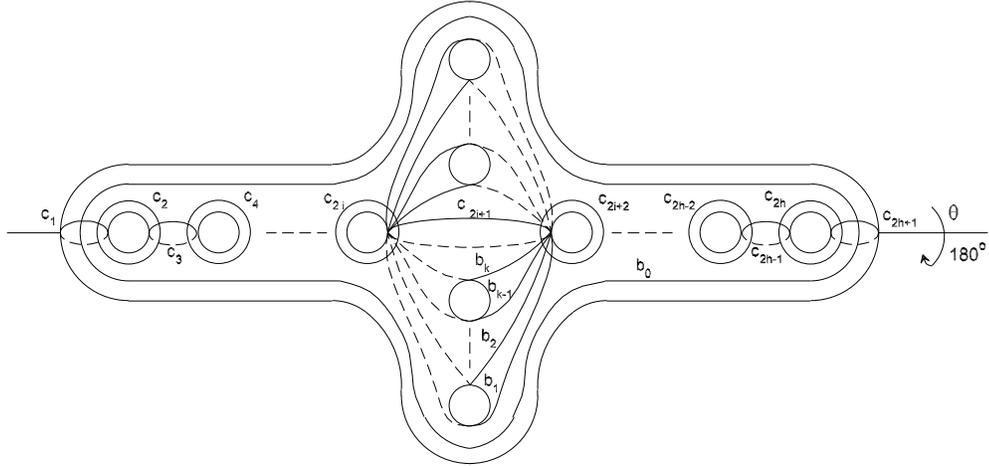,width=5.50in,clip=}
     \caption{The cycles used in the expression of $\theta$}
     \label{simplecase.fig}
 \end{figure}

Figure \ref{simplecase.fig} shows the cycles that are used in
expressing $\theta$ as a product of positive Dehn twists which is
stated in the next theorem.

\begin{thm} \label{simplecase.thm}
 The positive Dehn twist expression for the
involution $\theta$ shown in Figure \ref{simplecase.fig} is given
by
\[
\theta =c_{2i+2}\cdots c_{2h}c_{2h+1}c_{2i}\cdots
c_{2}c_{1}b_{0}c_{2h+1}c_{2h}\cdots c_{2i+2}c_{1}c_{2}\cdots
c_{2i}b_{1}b_{2}\cdots b_{k-1}b_{k}c_{2i+1}.
\]
\end{thm}

See \cite{Gu} for the proof.

\section{Main Theorem}

\subsection{The Involution $\theta$ on bounded surface}
\begin{figure}[htbp]
     \centering  \leavevmode
     \psfig{file=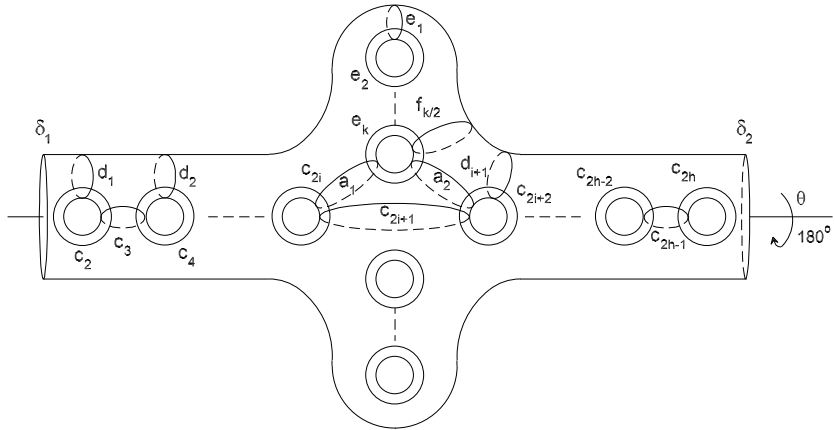,width=4.50in,clip=}
     \caption{The involution $\theta$ on bounded surface and a pants decomposition}
     \label{basecycleswithboundary.fig}
 \end{figure}

\begin{figure}[htbp]
     \centering  \leavevmode
     \psfig{file=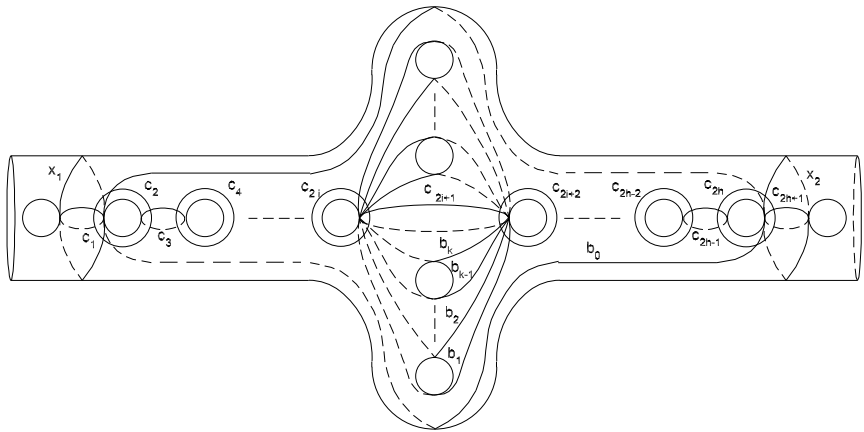,width=5.50in,clip=}
     \caption{The cycles realizing $\theta$ on the bounded surface in Figure \ref{basecycleswithboundary.fig} }
     \label{simplecasewithboundary.fig}
 \end{figure}

Consider the bounded surface $\Sigma_{h+k,2}$ in Figure
\ref{basecycleswithboundary.fig}, which is obtained from the
surface in Figure \ref{simplecase.fig} by removing a disk from
each end. Figure \ref{simplecasewithboundary.fig} is obtained from
Figure \ref{basecycleswithboundary.fig} by gluing a torus with two
boundary components on each end. The cycles shown in Figure
\ref{simplecasewithboundary.fig} realize the involution $\theta$
on the bounded surface in Figure \ref{basecycleswithboundary.fig}
as stated in proposition \ref{boundedcase.prop}.

The boundary components of the chosen pants decomposition shown in
Figure \ref{basecycleswithboundary.fig} will constitute the set of
cycles that will be mapped in order to prove proposition
\ref{boundedcase.prop}.

\begin{figure}[htbp]
     \centering  \leavevmode
     \psfig{file=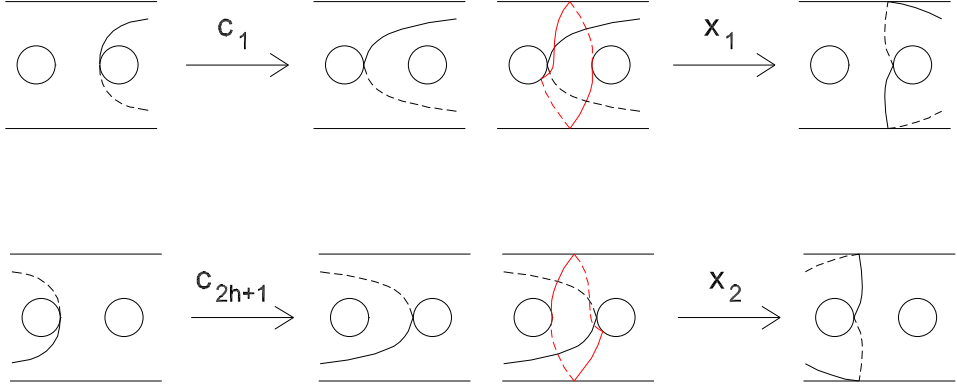,width=5.50in,clip=}
      \caption{The effect of $x_1c_1$ and
$x_2c_{2h+1}$ on the given arcs}
    \label{connectionlemma.fig}
 \end{figure}

Since the mapping of many of those boundary components will create
in the process cycles that will contain a piece of arc similar to
the two that are shown in the first column of Figure
\ref{connectionlemma.fig}, we will show the mappings of these
segments separately once and use their images in the last column
of the same figure to avoid repetition, whenever necessary in the
proof of the proposition.

Each row in the following lemma shows the mapping of one of the
two types of segments that will occur several times in the proof
of proposition \ref{boundedcase.prop} as mentioned above.

\begin{lem} \label{connectionlemma.lem} The action of the Dehn twists $x_1c_1$ and
$x_2c_{2h+1}$ on the arcs shown in the first column of Figure
\ref{connectionlemma.fig} are as shown in the last column of the
same figure.
\end{lem}

\begin{prop} \label{boundedcase.prop} The positive Dehn twist expression
for the involution $\theta$ defined on the bounded surface
$\Sigma_{h+k,2}$ shown in Figure \ref{basecycleswithboundary.fig}
is given by
\[
\theta =c_{2i+2}\cdots c_{2h}c_{2i}\cdots
c_{2}x_{2}c_{2h+1}x_{1}c_{1}b_{0}c_{2h}\cdots c_{2i+2}c_{2}\cdots
c_{2i}b_{1}b_{2}\cdots b_{k-1}b_{k}c_{2i+1},
\]
where the cycles in the expression are as shown in Figure
\ref{simplecasewithboundary.fig}.
\end{prop}

\noindent {\bf Proof:} Figure \ref{basecycleswithboundary.fig}
shows a pants decomposition for the bounded surface on which
$\theta$ is defined. We will show that the given Dehn twist
expression in the proposition maps the boundary components of each
pair of pants to their images under $\theta$. This will guarantee
the mapping of the interior points of each pair of pants
accordingly, due to the fact that each twist in the expression is
a homeomorphism of the surface onto itself.

The same idea was used in proving theorem \ref{simplecase.thm} in
\cite{Gu} for the closed surface $\Sigma_{h+k}$ and the mapping of
each boundary cycle was shown there in detail, up to symmetry.
Even though the surface subject to this proposition is not closed,
there are several figures that are identical for both cases.
Therefore, for a given boundary component, instead of repeating
verbatim copy of the figures in its mapping from \cite{Gu}, we
will skip a few from the beginning and continue from where the
different cycles begin to appear. The reader is referred to that
article for the details of the mappings that are skipped here.

The boundary components of the chosen pants decomposition in
Figure \ref{basecycleswithboundary.fig} can be summarized as $c_i,
i=1,\ldots,2h+1, d_i,i=2,\ldots,h-1, e_i, i=1,\ldots,2k+1,
f_i,i=2,\ldots,k-1, a_1,a_2,\delta_1$ and $\delta_2$ along with
some additional cycles.

\begin{figure}[htbp]
     \centering  \leavevmode
     \psfig{file=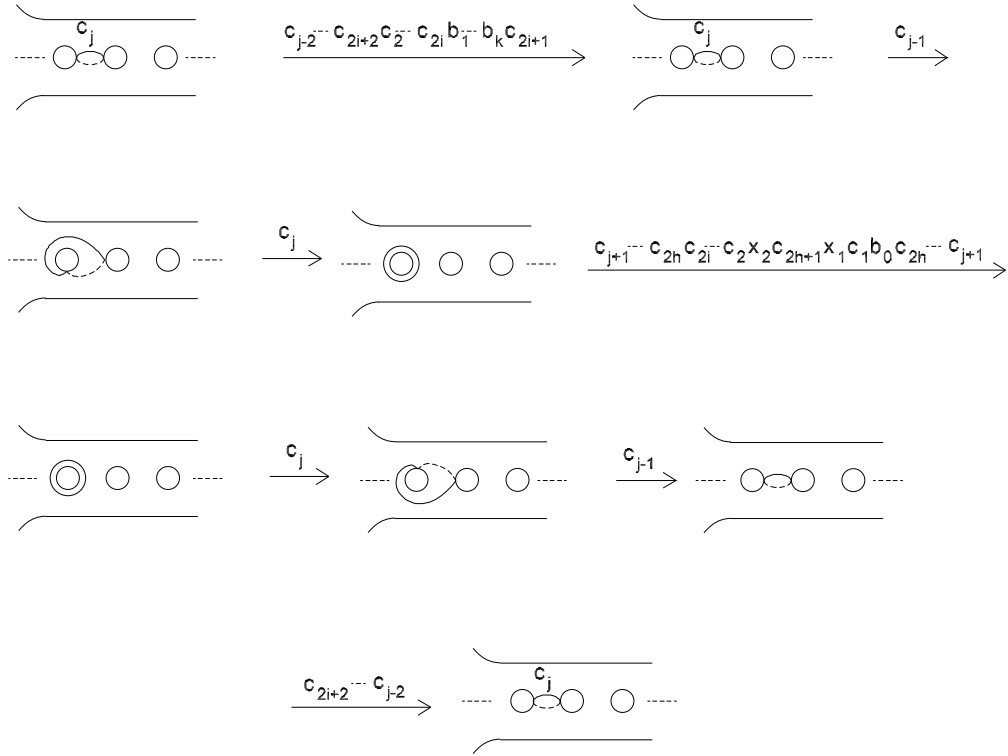,width=8.0in,clip=}
     \caption{The mapping of $c_j$}
     \label{mappingofcjwb.fig}
 \end{figure}

We will begin with the mapping of $c_j$ for $j-$ odd and $2i+3\leq
j< 2h$, Figure \ref{mappingofcjwb.fig}. The proof for $j-$ even,
including $j=2$ and $j=2h$, is similar and was shown in \cite{Gu}.
The mappings of $c_{2i+1}$ and $c_{2i+2}$ will  be shown
separately.

The long expressions in Figure \ref{mappingofcjwb.fig} are due the
fact that all the twists they contain miss the cycle that appears
in the previous step. The figure shows all the steps there are.

\begin{figure}[htbp]
     \centering  \leavevmode
     \psfig{file=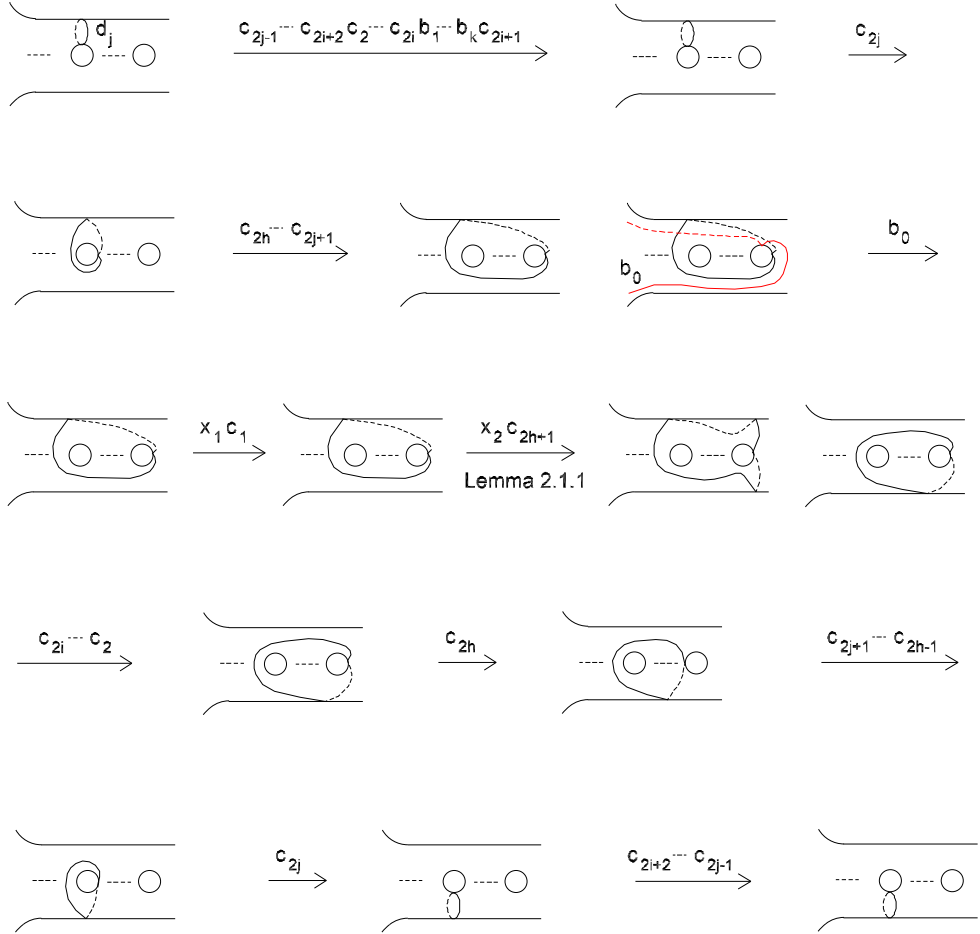,width=10.0in,clip=}
      \caption{The mapping of $d_j$}
    \label{mappingofdjwb.fig}
 \end{figure}

Wee see the mapping of $d_j$ in Figure \ref{mappingofdjwb.fig},
which is the same for $j=i+1,\ldots,h$. The twist about $b_0$
leaves the curve it is applied to unchanged because their
intersection number is 0 as seen in the end of the second line.
The result of application of the twists $x_2c_{2h+1}$ is obtained
according to Lemma \ref{connectionlemma.lem}, therefore only the
right end portion of the cycle to which they are applied is
modified in the third line. The cycle in the end of the third line
 is isotopic to the previous one because it is obtained simply by
retracting the portion that falls under the surface. The mapping
of $d_j$ for $j=1,\ldots,i$ is similar due to symmetry and is
omitted.

\begin{figure}[htbp]
     \centering  \leavevmode
     \psfig{file=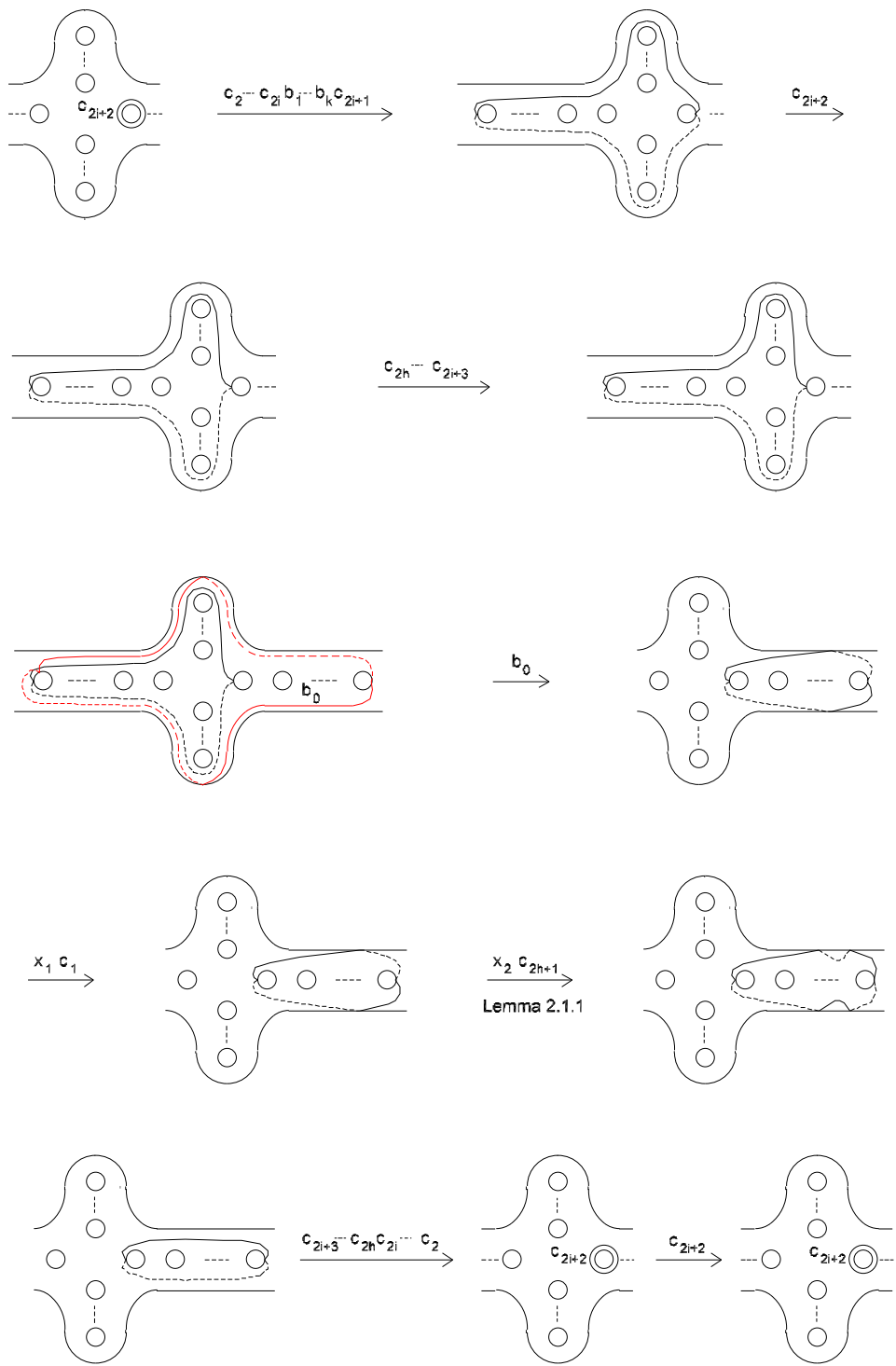,width=14.50in,clip=}
      \caption{The mapping of $c_{2i+2}$}
    \label{c2iplus2wb.fig}
 \end{figure}

Figure \ref{c2iplus2wb.fig} shows the mapping of $ c_{2i+2}$. The
details of the applications of the twists $c_2\cdots
c_{2i}b_{1}\cdots b_{k}c_{2i+1}$ in the first line are skipped and
can be found in \cite{Gu}. Note the use of Lemma
\ref{connectionlemma.lem} in the second line from the bottom. The
first cycle of the last line is isotopic to the one that appears
just before.

\clearpage

\begin{figure}[htbp]
     \centering  \leavevmode
     \psfig{file=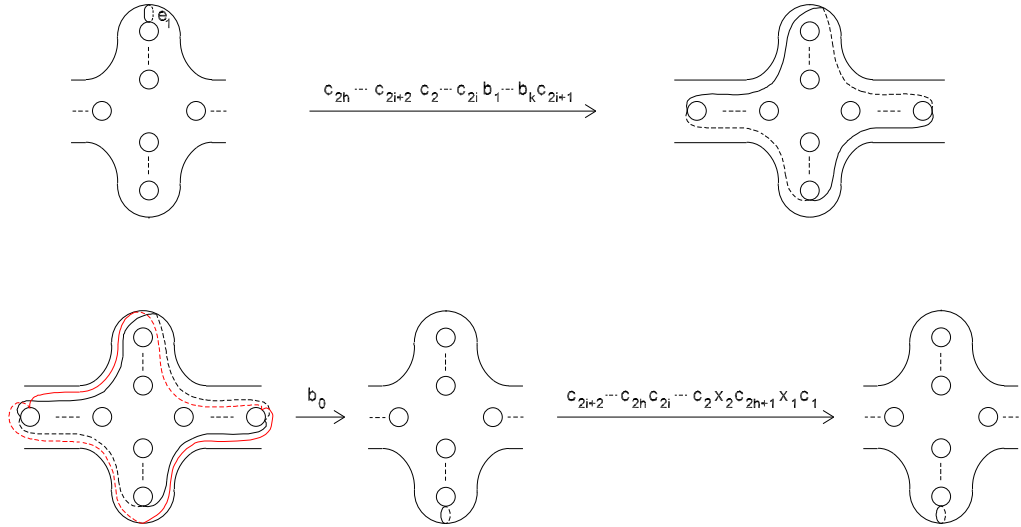,width=6.0in,clip=}
      \caption{The mapping of $e_1$}
    \label{mappingofe1wb.fig}
 \end{figure}

\begin{figure}[htbp]
     \centering  \leavevmode
     \psfig{file=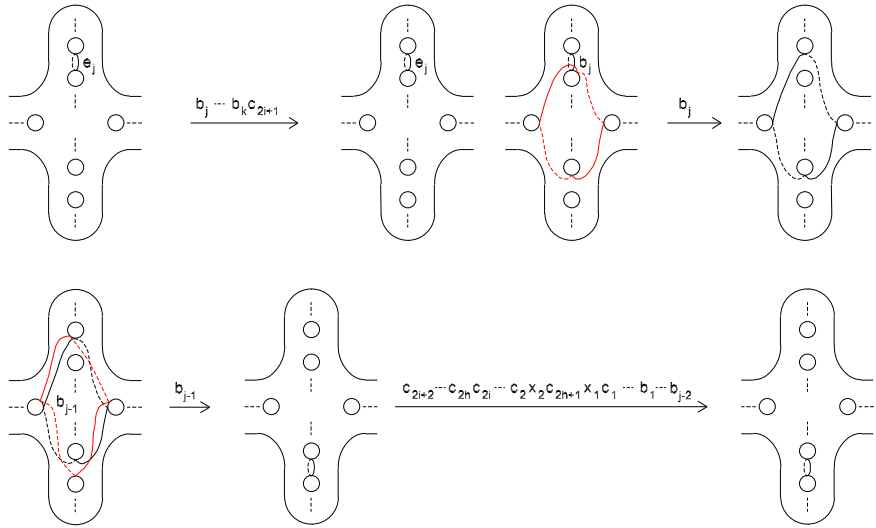,width=5.50in,clip=}
      \caption{The mapping of $e_j$}
    \label{mappingofejwb.fig}
 \end{figure}

The only curves that are effective in the mappings of $e_j$ are
$b_j$ and $b_{j-1},j=1,\ldots,k$. Figures \ref{mappingofe1wb.fig}
and \ref{mappingofejwb.fig} show the mapping of $e_j$ for $j-$
odd. The twists in the long expressions all miss the curves that
come before them.

\begin{figure}[htbp]
     \centering  \leavevmode
     \psfig{file=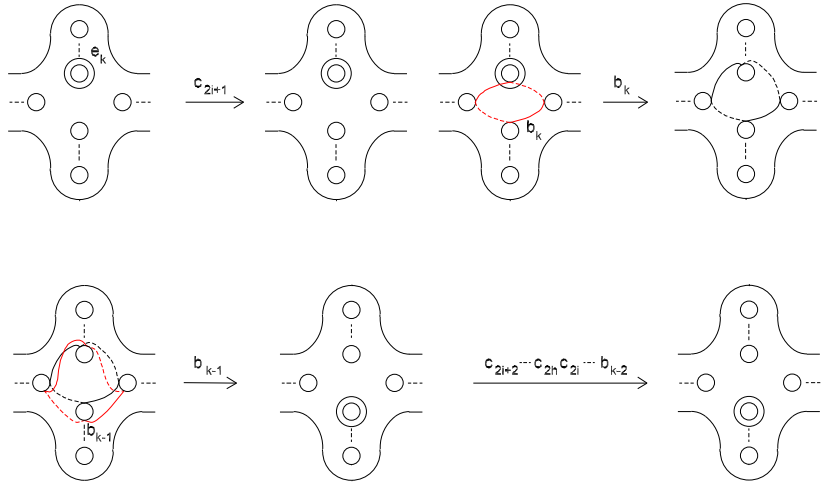,width=5.50in,clip=}
      \caption{The mapping of $e_k$}
    \label{mappingofekwb.fig}
 \end{figure}

The mapping of $e_k$ is a typical example for the mapping of
$e_j,j-$ even, which is shown in Figure \ref{mappingofekwb.fig}.

\begin{figure}[htbp]
     \centering  \leavevmode
     \psfig{file=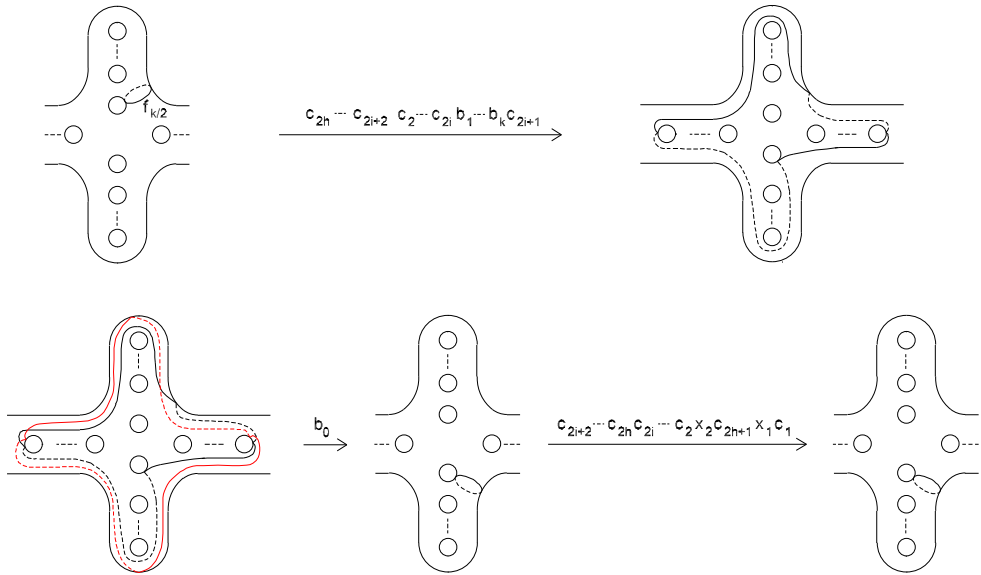,width=6.0in,clip=}
      \caption{The mapping of $f_j$}
    \label{mappingoffjwb.fig}
 \end{figure}

The mapping of $f_{k/2}$ is shown in Figure
\ref{mappingoffjwb.fig}. The details of the applications of the
twists $c_{2h}\cdots c_{2i+2}c_2\cdots c_{2i}b_{1}\cdots
b_{k}c_{2i+1}$ in the first line are skipped and can be found in
\cite{Gu}. The mappings of $f_j$ for $j=2,\ldots,k/2-1$ are
similar to that of $f_{k/2}$. Note that $f_1$ is the same as
$e_1$.

\begin{figure}[htbp]
     \centering  \leavevmode
     \psfig{file=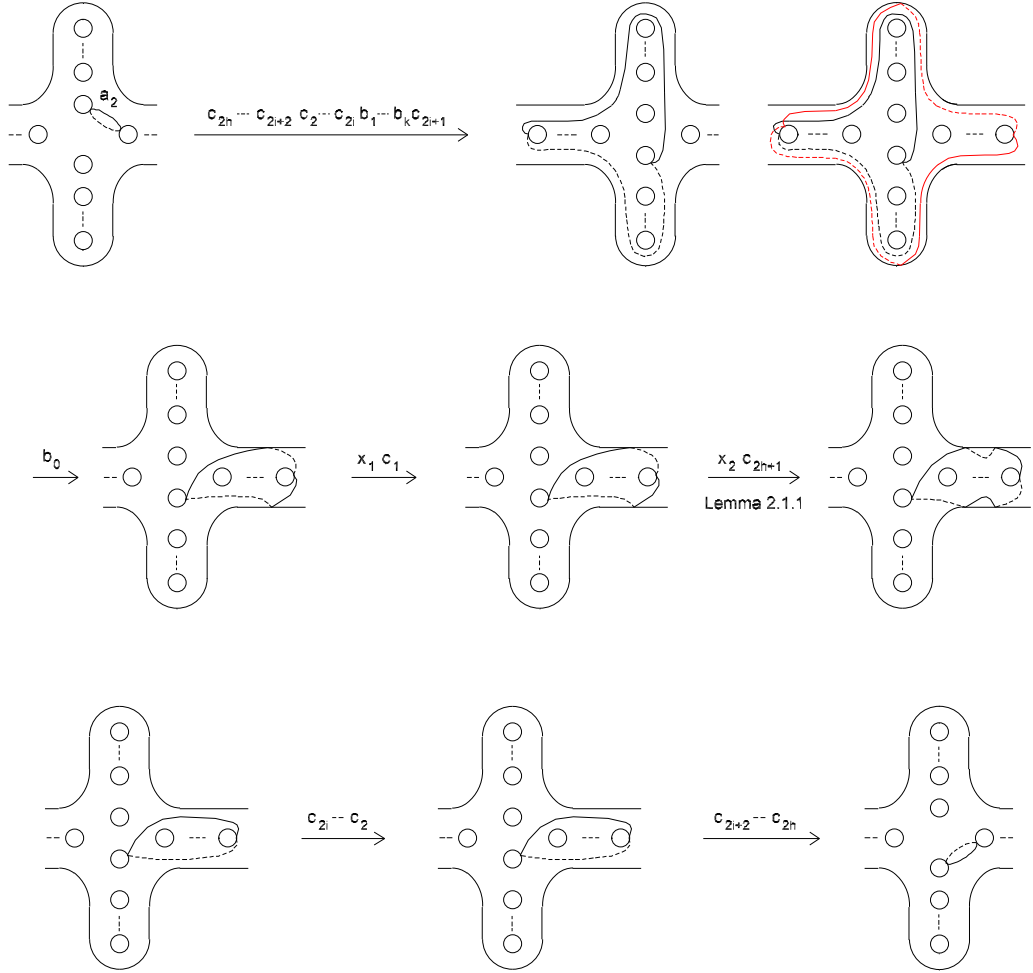,width=9.50in,clip=}
      \caption{The mapping of $a_2$}
    \label{mappingofa2wb.fig}
 \end{figure}

Figure \ref{mappingofa2wb.fig} shows the mapping $a_2$ and the
mapping of $a_1$ is symmetrical to it. In this figure the details
of the applications of the twists $c_{2h}\cdots c_{2i+2}c_2\cdots
c_{2i}b_{1}\cdots b_{k}c_{2i+1}$ are skipped also in the first
line. Lemma \ref{connectionlemma.lem} is used in the second line
and the resulting curve from that is isotopic to the curve in the
beginning of the third line.

\begin{figure}[htbp]
     \centering  \leavevmode
     \psfig{file=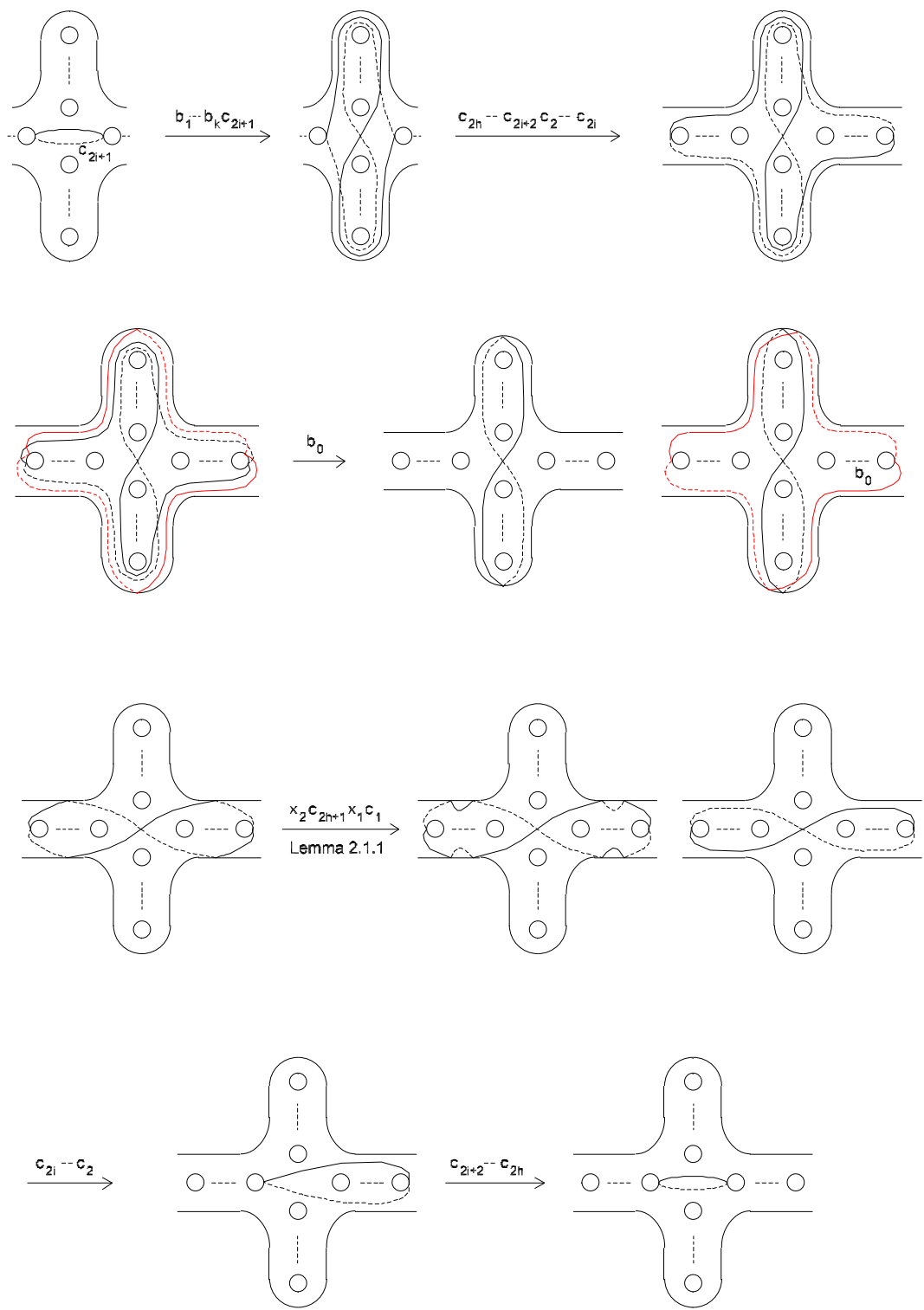,width=15.0in,clip=}
      \caption{The mapping of $c_{2i+1}$}
    \label{c2iplus1wb.fig}
 \end{figure}

In Figure \ref{c2iplus1wb.fig} we see the mapping of $c_{2i+1}$.
In this figure, too, the details of the applications of
$b_{1}\cdots b_{k}c_{2i+1}$ and $c_{2h}\cdots c_{2i+2}c_2\cdots
c_{2i}$ are skipped in the first line. Lemma
\ref{connectionlemma.lem} is used twice in the third line and the
last figure in that line is isotopic to the one that is resulting
from the application of the lemma. Note that $b_0$ has
intersection number 2 with the curve it is applied to; therefore,
the result of the twist about $b_0$ is found by taking their
product twice.

 Finally, we will show the mapping of $\delta_1$, which is
 essentially the same as that of $\delta_2$ due to symmetry.

\begin{figure}[htbp]
     \centering  \leavevmode
     \psfig{file=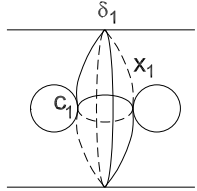,width=2.0in,clip=}
      \caption{ $\delta_1$ and the cycles that are effective in its mapping}
    \label{delta1}
 \end{figure}

\begin{figure}[htbp]
     \centering  \leavevmode
     \psfig{file=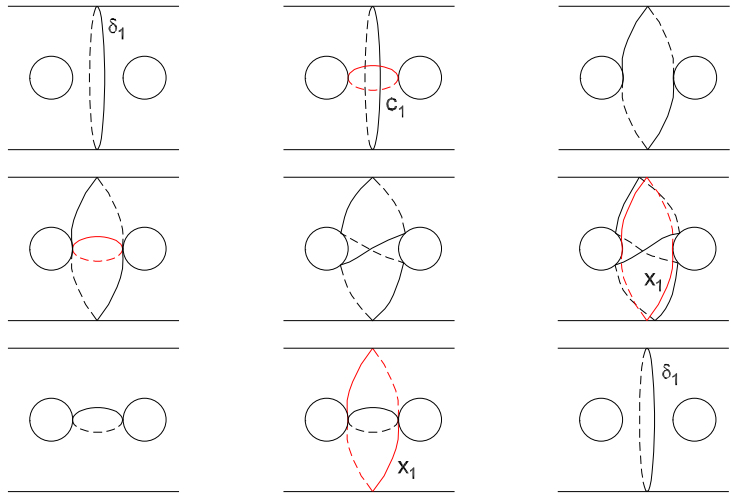,width=5.0in,clip=}
      \caption{The mapping of $\delta_1$}
    \label{mappingofdelta1.fig}
 \end{figure}

The only cycles that take part in the mapping of $\delta_1$
 are $c_1$ and $x_1$, as shown in Figure \ref{delta1}. All the cycles that
 come before $c_1$ miss $\delta_1$ as well as the ones that come
 after $x_1$ and  $x_1c_1$ fixes $\delta_1$ point-wise. This
 is shown in Figure \ref{mappingofdelta1.fig}.

The intersection number of $c_1$ and
 $\delta_1$ is 2; therefore $c_1(\delta_1)=c_1^2\delta_1$,
 namely the product of $c_1$ and  $\delta_1$ twice.
$c_1(\delta_1)$  is the second cycle in the second row of Figure
\ref{mappingofdelta1.fig}.
 The intersection number of  $c_1(\delta_1)$ and $x_1$ is also 2;
 therefore $x_1(t_1(\delta_1))=x_1^2t_1(\delta_1)$,
 which is $\delta_1$ as seen in the last row of the same figure. Therefore
 $x_1c_1(\delta_1)=\delta_1$, as claimed.

\begin{figure}[htbp]
     \centering  \leavevmode
     \psfig{file=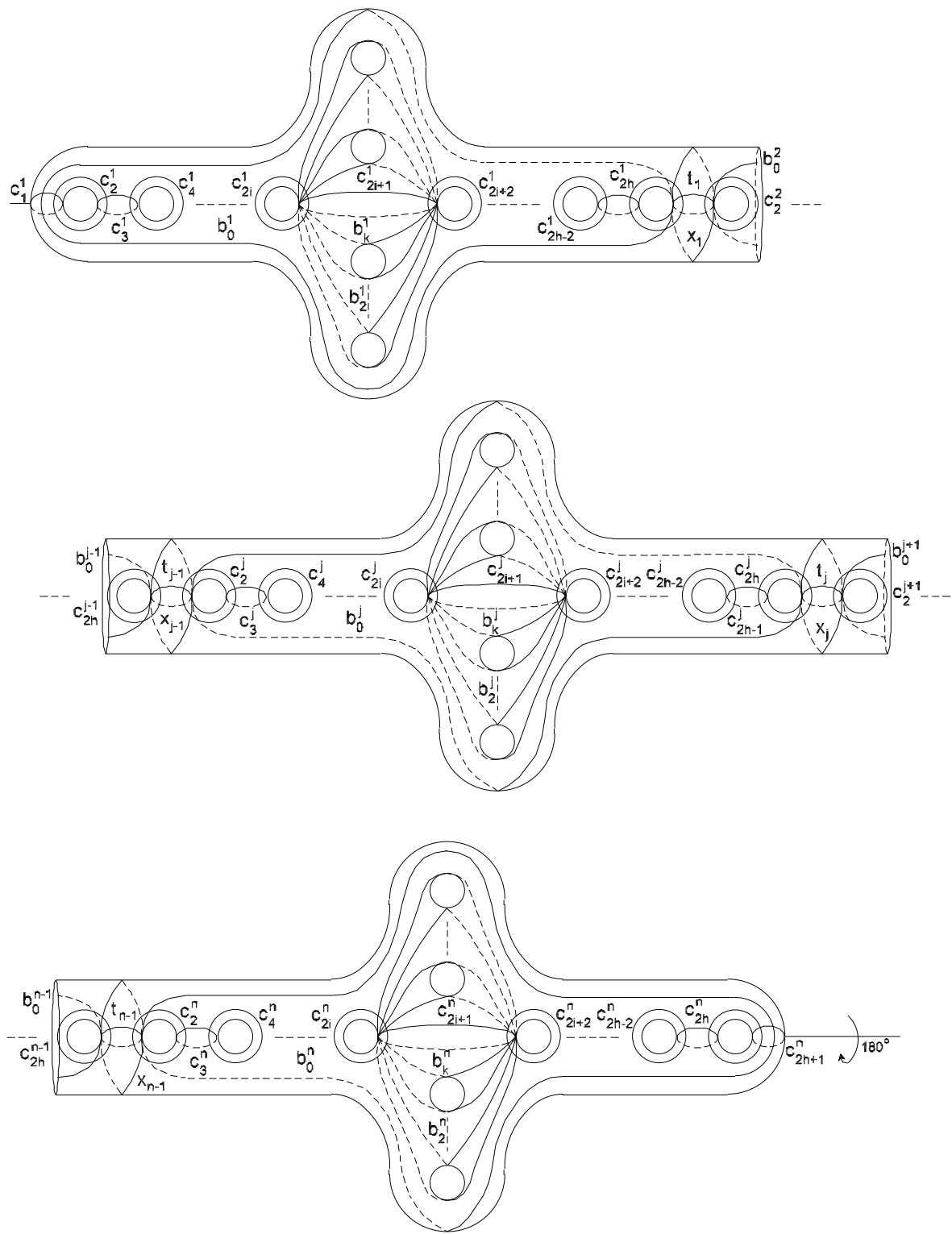,width=15.0in,clip=}
      \caption{General case}
    \label{multiple.fig}
 \end{figure}

This concludes the proof of Proposition \ref{boundedcase.prop}.
Now we can prove the main theorem, which is the generalization of
Proposition \ref{boundedcase.prop} to a surface that is obtained
by gluing $n$ copies of bounded surfaces as in Figure
\ref{basecycleswithboundary.fig} together along four-holed spheres
in a sequence. Each copy in that sequence will then have two
boundary components except for the first and the last copies,
which will have only one boundary component each as shown in
Figure \ref{multiple.fig}.

\clearpage

In order to simplify the Dehn twist expression for the general
case it will be necessary to group the twists in each copy and
give suggestive names to them. We will also pay attention to the
direction in which the horizontal twists are progressing.

The label of each twist group will carry the following
information: Which copy the twist group is in (upper index),
whether the group is on the right-hand side or on the left-hand
side or in the middle section of the respective copy (name of the
group), the direction in which the twists are multiplied when the
group consists of horizontal twists (lower index). The twists
along the cycle $b_0$ will not be included in any group. The
following is the list of the identifications, except for the first
and the last copies:

 \begin{eqnarray*}
 l^j_i&=&c^j_{2i}\cdots c^j_2,\\ l^j_o&=&c^j_{2}\cdots
c^j_{2i},\\ r^j_i&=&c^j_{2i+2}\cdots c^j_{2h},
\\ r^j_o&=&c^j_{2h}\cdots c^j_{2i+2}, \\ m^j&=&b^j_1\cdots
b^j_{k}c^j_{2i+1}.
\end{eqnarray*}

The first two lines would be different in the first copy and the
third and the fourth lines would be different in the last copy:

 \begin{eqnarray*}
 l^1_i&=&c^1_{2i}\cdots c^1_1,\\ l^1_o&=&c^1_{1}\cdots
c^1_{2i},\\ r^n_i&=&c^n_{2i+2}\cdots c^n_{2h+1},
\\ r^n_o&=&c^n_{2h+1}\cdots c^n_{2i+2}.
\end{eqnarray*}

 Basically $l^j_i$ is the \emph{inward} product of the horizontal
twists on the \emph{left} hand side of the $j^{th}$ copy, namely
their product taken towards the center of the $j^{th}$ copy.
Similarly $l^j_o$ is the \emph{outward} product of the twists on
the \emph{left} hand side of the $j^{th}$ copy, namely their
product taken away from the center of the $j^{th}$ copy. The
definitions of $r^j_i$ and $r^j_o$ use the same idea. $m^j$
represents the product of the twists in the \emph{middle} section
of the $j^{th}$ copy as it appears in Theorem \ref{simplecase.thm}

The twists $c^j_{2i},c^1_{2i},c^j_{2i+2}$ and  $c^n_{2i+2}$ should
actually be written as $c^j_{2i_j},c^1_{2i_1},c^j_{2i_j+2}$ and
$c^n_{2i_n+2}$ as the subindex $i$ will be different for each copy
but we are not showing this dependence on $j$ to keep the notation
simple.

 Using the notation described above we can write the positive Dehn
 twist product for the involution shown in Figure
 \ref{multiple.fig}:

\begin{thm} The positive Dehn twist product for the involution
\label{main.thm} $\theta$ shown in Figure \ref{multiple.fig} is
\[r^n_il^n_ir^{n-1}_il^{n-1}_ix_{n-1}t_{n-1}b^n_0r^n_ol^n_om^n\cdots
m^4r^2_il^2_ix_2t_2b^3_0r^3_ol^3_om^3r^1_il^1_ix_1t_1b^2_0r^2_ol^2_om^2b^1_0r^1_ol^1_om^1.
\]
\end{thm}

To reduce the notation in the above expression for $\theta$
further let's let
 \begin{eqnarray*}
   Y^j_o&=&r^j_ol^j_om^j,\\
   Y^j_i&=&r^j_il^j_i,\\
   X_j&=&x_jt_j.
 \end{eqnarray*}

 Then $\theta$ can be rewritten as
 \[Y^n_iY^{n-1}_iX_{n-1}b^n_0Y^n_o\cdots
 b^4_0Y^4_oY^2_iX_2b^3_0Y^3_oY^1_iX_1b^2_0Y^2_ob^1_0Y^1_o.
\]

Using the product notation we obtain:
\[\theta =Y^n_i\prod_{j=2}^{n}\left(Y^{j-1}_iX_{j-1}b^j_0Y^j_o\right)b^1_0Y^1_o
\]

The product sign here will mean multiplication from right to left,
contrary to its usual meaning, in agreement with the earlier
expressions.\\

\begin{figure}[htbp]
     \centering  \leavevmode
     \psfig{file=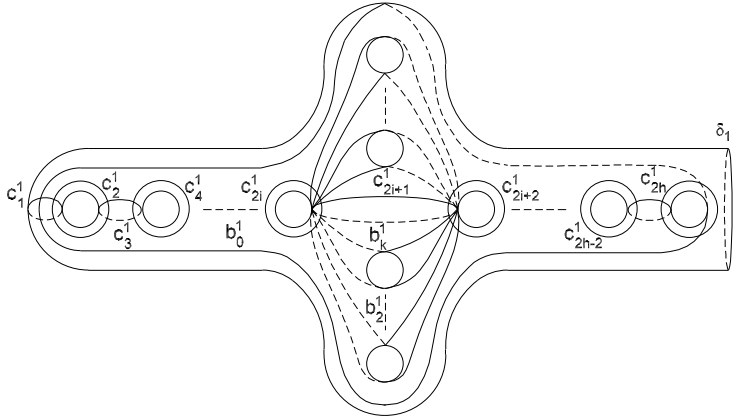,width=4.50in,clip=}
      \caption{Initial step}
    \label{nequals1.fig}
 \end{figure}

\noindent {\bf Proof:} The proof is by induction. To show the
effect of $\theta$ on the first copy we set $n=2$ in the product
sign and get the expression
\[
Y^{1}_iX_{1}b^2_0Y^2_ob^1_0Y^1_o,
\]
 which is equal to

\[
Y^{1}_iX_{1}b^1_0Y^1_ob^2_0Y^2_o,
\]
because none of the cycles in $b^1_0Y^1_o$ intersects any cycle in
$b^2_0Y^2_o$. Therefore what we have is
\[r^1_il^1_ix_1t_1b^1_0r^1_ol^1_om^1b^2_0r^2_ol^2_om^2,
\]
which can be reduced to

\[r^1_il^1_ix_1t_1b^1_0r^1_ol^1_om^1
\]
because the expression $b^2_0r^2_ol^2_om^2$ has no effect on the
bounded surface in Figure \ref{nequals1.fig}.

The explicit version of $r^1_il^1_ix_1t_1b^1_0r^1_ol^1_om^1$ is
\[
c_{2i+2}^1\cdots c_{2h}^1c_{2i}^1\cdots
c_{1}^1x_{1}^1t_{1}^1b_{0}^1c_{2h}^1\cdots
c_{2i+2}^1c_1^1c_{2}^1\cdots c_{2i}^1b_{1}^1b_{2}^1\cdots
b_{k-1}^1b_{k}^1c_{2i+1}^1,
\]

This is a special case of the expression in Proposition
\ref{boundedcase.prop} for the surface with one boundary
component. Therefore the effect of the above expression on the
bounded surface in Figure \ref{nequals1.fig} is that of $\theta$
in Proposition \ref{boundedcase.prop}.

 Now suppose that
\[\prod_{u=2}^{j}\left(Y^{u-1}_iX_{u-1}b^u_0Y^u_o\right)b^1_0Y^1_o
\]
realizes the involution $\theta$ on the first $j-1$ copies of the
surface in Figure \ref{multiple.fig}. Consider now
\[\prod_{u=2}^{j+1}\left(Y^{u-1}_iX_{u-1}b^u_0Y^u_o\right)b^1_0Y^1_o,
\]
which is equal to
\[Y^{j}_iX_{j}b^{j+1}_0Y^{j+1}_o\prod_{u=2}^{j}\left(Y^{u-1}_iX_{u-1}b^u_0Y^u_o\right)b^1_0Y^1_o.
\]

\begin{figure}[htbp]
     \centering  \leavevmode
     \psfig{file=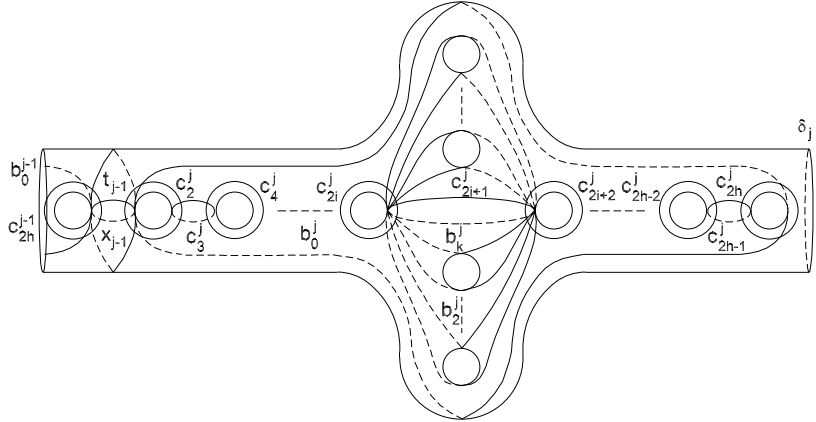,width=4.50in,clip=}
      \caption{Inductive step}
    \label{nequalsj.fig}
 \end{figure}

The first observation we have to make is, the expression
$Y^{j}_iX_{j}b^{j+1}_0Y^{j+1}_o$ leaves the first $j-1$ copies
with boundary $\delta_{j-1}$ unaltered, because all of the twists
it contains are about cycles that lie completely to the right of
$\delta_{j-1}$. Now, in order to see the effect of the inductive
step on the surface in Figure \ref{nequalsj.fig} let's release the
last term in the product to get
\[Y^{j}_iX_{j}b^{j+1}_0Y^{j+1}_oY^{j-1}_iX_{j-1}b^{j}_0Y^{j}_o
\prod_{u=2}^{j-1}\left(Y^{u-1}_iX_{u-1}b^u_0Y^u_o\right)b^1_0Y^1_o.
\]

The part of this expression that will be effective in the mapping
of the $j^{th}$ copy is contained in
$Y^{j}_iX_{j}b^{j+1}_0Y^{j+1}_oY^{j-1}_iX_{j-1}b^{j}_0Y^{j}_o$.

Using the commutativity relation between the terms that do not
intersect we can rewrite this as
$Y^{j-1}_iY^{j}_iX_{j}X_{j-1}b^{j}_0Y^{j}_ob^{j+1}_0Y^{j+1}_o$,
just to bring the terms that we need together. To be precise, the
twists contained in $Y^{j}_iX_{j}X_{j-1}b^{j}_0Y^{j}_o$ are the
ones that will realize the effect of $\theta$ on the $j^{th}$
copy. Writing them explicitly, we get
\[
c^j_{2i+2}\cdots c^j_{2h}c^j_{2i}\cdots
c^j_{2}x_{j}t_{j}x_{j-1}t_{j-1}b^j_{0}c^j_{2h}\cdots
c^j_{2i+2}c^j_{2}\cdots c^j_{2i}b^j_{1}b^j_{2}\cdots
b^j_{k-1}b^j_{k}c^j_{2i+1},
\]

which is exactly the expression in Proposition
\ref{boundedcase.prop} adapted for the $j^{th}$ copy, with the
identifications $c_1=t_{j-1},c_{2h+1}=t_j,x_1=x_{j-1},x_2=x_j$.
This proves the inductive step.

To complete the proof we need to point out to the mapping of the
last copy. Recall the expression for $\theta$
\[\theta
=Y^n_i\prod_{j=2}^{n}\left(Y^{j-1}_iX_{j-1}b^j_0Y^j_o\right)b^1_0Y^1_o.\]
Releasing the last term in the product sign we get
\[
Y^n_iY^{n-1}_iX_{n-1}b^n_0Y^n_o\prod_{j=2}^{n-1}\left(Y^{j-1}_iX_{j-1}b^j_0Y^j_o\right)b^1_0Y^1_o
\]
Since $Y^{n-1}_i$ has no effect on the $n^{th}$ copy we have only
$Y^n_iX_{n-1}b^n_0Y^n_o$ realizing $\theta$ on the last copy.
Writing them explicitly we get
\[
c^n_{2i+2}\cdots c^n_{2h}c^n_{2h+1}c^n_{2i}\cdots
c^n_{2}x_{n-1}t_{n-1}b^n_{0}c^n_{2h+1}c^n_{2h}\cdots
c^n_{2i+2}c^n_{2}\cdots c^n_{2i}b^n_{1}b^n_{2}\cdots
b^n_{k-1}b^n_{k}c^n_{2i+1},
\]

\begin{figure}[htbp]
     \centering  \leavevmode
     \psfig{file=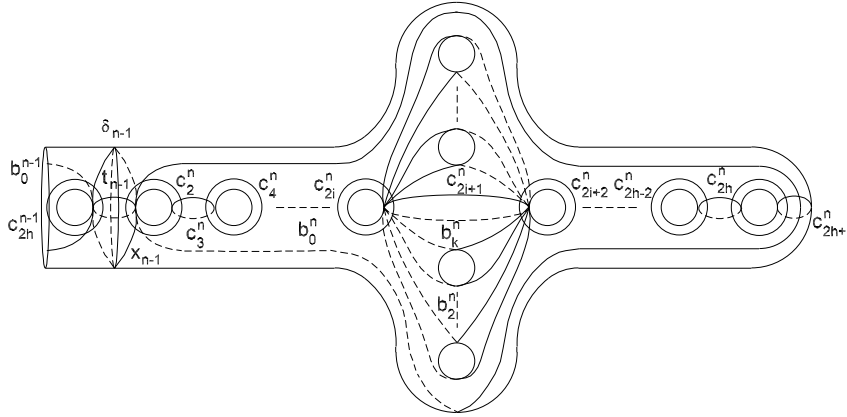,width=4.50in,clip=}
      \caption{Last copy}
    \label{nequalsn.fig}
 \end{figure}

This is, again, a special case of the formula in Proposition
\ref{boundedcase.prop} adapted for the surface with one boundary
component seen in Figure \ref{nequalsn.fig}.

Although it is not needed, we will also include the mapping of the
cycle $t_j$ in the proof.

\begin{figure}[htbp]
     \centering  \leavevmode
     \psfig{file=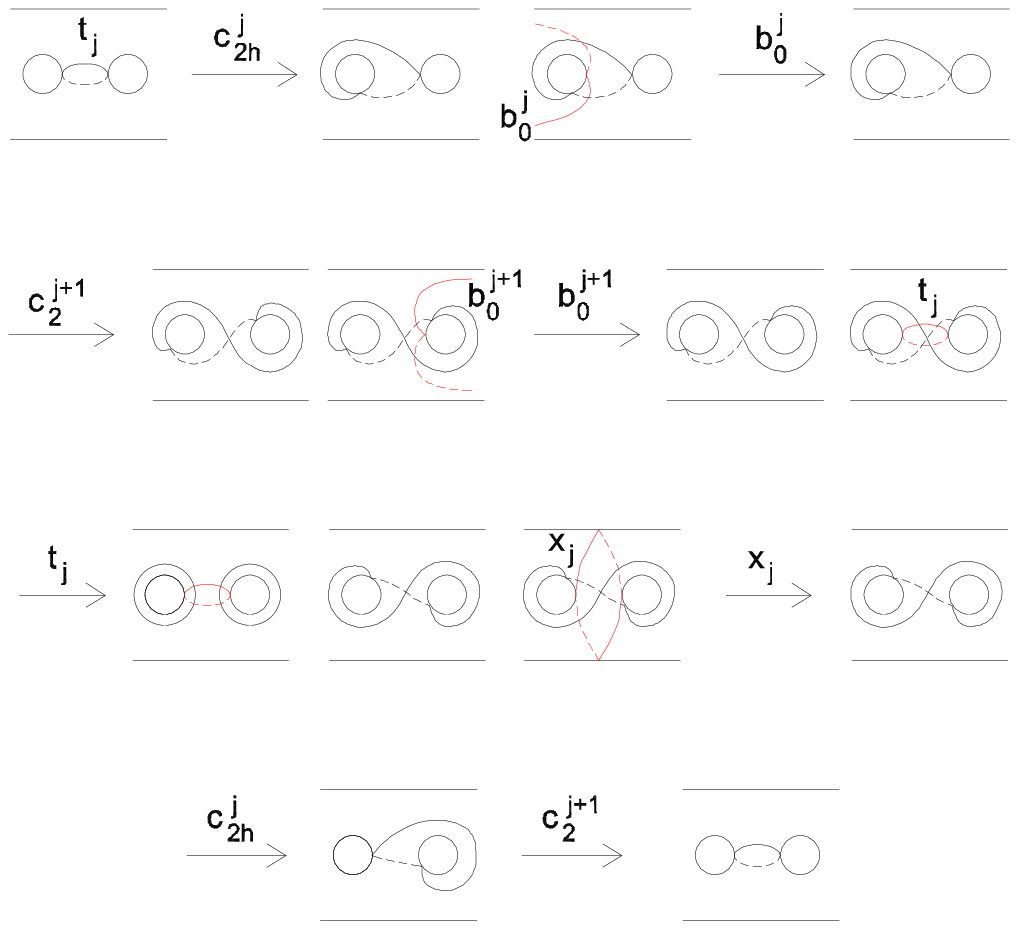,width=8.50in,clip=}
      \caption{The mapping of $t_j$}
    \label{mappingoftjwb.fig}
 \end{figure}

Figure \ref{mappingoftjwb.fig} shows the mapping of the curve
$t_j,j=1\ldots,n-1$. To understand the steps in that figure let's
write
\[
\theta=\cdots
Y^{j+1}_iX_{j+1}b^{j+2}_0Y^{j+2}_oY^{j}_iX_{j}b^{j+1}_0Y^{j+1}_oY^{j-1}_iX_{j-1}b^j_0Y^j_o\cdots.
\]
The first term in the expression for $\theta$ that will not miss
$t_j$ is $Y^j_o=r^j_ol^j_om^j$. In fact all the twists in $Y^j_o$
will miss $t_j$ except for the last twist in $r^j_o,$ which is
$c^j_{2h}$. The next twist is $b_0^j$ and it will leave the result
of the previous twist unaltered as shown in the first line of
Figure \ref{mappingoftjwb.fig}. So does the expression
$Y^{j-1}_iX_{j-1}$ because all the twists they contain miss the
same result. The effect of the next term
$Y^{j+1}_o=r^{j+1}_ol^{j+1}_om^{j+1}$ on the current cycle is
performed by the twist $c^{j+1}_2$ which is contained in
$l^{j+1}_o$. The result from this twist is seen in the second line
of the figure. This movement causes the next twist to miss the
current result, namely $b^{j+1}_0$ leaves it unaltered as shown in
the second line. The first twist $t_j$ in $X_j=x_jt_j$ has two
intersection points with the current cycle and the result from its
application is seen in the first half of the third line. The cycle
$x_j$ doesn't intersect the result from twisting about $t_j$,
therefore it has no effect on it as indicated in the end of the
third line. The following term $Y^j_i=r^j_il^j_i$ has only one
cycle that will intersect the cycle that is missed by $x_j$ in the
third line, i.e., $c^j_{2h}$ that lies in $r^j_i$. Its effect on
the current cycle is seen in the beginning of the last line. All
the twists contained in the sequence of terms
$X_{j+1}b^{j+2}_0Y^{j+2}_o$ following $Y^{j}_i$ miss the first
cycle in the last line. The next cycle that will not miss it is
$c_2^{j+1}$, which is contained in $l^{j+1}_i$ of
$Y^{j+1}_i=r^{j+1}_il^{j+1}_i$. The rest of the twists miss the
last cycle in Figure \ref{mappingoftjwb.fig}, therefore $t_j$ is
fixed point-wise under the action of the expression for $\theta$,
as expected.

\begin{cor} \label{hyperelliptic.cor1} Let $\theta$ be expressed as in Theorem \ref{simplecase.thm}.
By setting $k=0$ we obtain the positive Dehn twist expression
\[
i =c_{2i+2}\cdots c_{2h}c_{2h+1}c_{2i}\cdots
c_{2}c_{1}b_{0}c_{2h+1}c_{2h}\cdots c_{2i+2}c_{1}c_{2}\cdots
c_{2i}c_{2i+1}
\]
for the hyperelliptic involution.

\end{cor}

\begin{figure}[htbp]
     \centering  \leavevmode
     \psfig{file=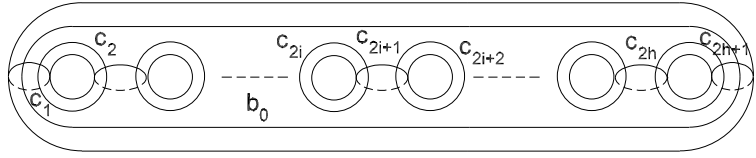,width=4.0in,clip=}
      \caption{A positive Dehn twist expression for the hyperelliptic involution}
    \label{corollary1.fig}
 \end{figure}

\noindent {\bf Proof:} We will give an algebraic proof for this
fact. First, observe that $b_0=c_1\cdots c_{2h}(c_{2h+1})$. Here
we abuse the notation and use the same notation for the cycles and
the twists. $b_0=c_1\cdots c_{2h}(c_{2h+1})$ means that the
sequence of twists $c_1\cdots c_{2h}$ are applied to the cycle
$c_{2h+1}$ and the cycle $b_0$ is obtained as the result of that.
Therefore the twist about $b_0$ is obtained from the twist about
$c_{2h+1}$ by conjugation by $c_1\cdots c_{2h}$ by a well-known
fact \cite{Gu}, i.e.,
\[b_0=c_1\cdots c_{2h}c_{2h+1} (c_1\cdots c_{2h})^{-1}.
\]

Substituting this in the expression for $i$ stated in Corollary
\ref{hyperelliptic.cor1} we obtain
\[
c_{2i+2}\cdots c_{2h}c_{2h+1}c_{2i}\cdots c_{2}c_{1}c_1\cdots
c_{2h}c_{2h+1} (c_1\cdots c_{2h})^{-1}c_{2h+1}c_{2h}\cdots
c_{2i+2}c_{1}c_{2}\cdots c_{2i}c_{2i+1},
\]
\[
c_{2i+2}\cdots c_{2h}c_{2h+1}c_{2i}\cdots c_{2}c_{1}c_1\cdots
c_{2h}c_{2h+1} c_{2h}^{-1}\cdots c_1^{-1} c_{2h+1}c_{2h}\cdots
c_{2i+2}c_{1}c_{2}\cdots c_{2i}c_{2i+1}.
\]

Recall that $c_i$ and $c_j$ commute if $|i-j|>1$. Using this we
can write
\[c_{2h}^{-1}\cdots c_1^{-1} c_{2h+1}c_{2h}\cdots
c_{2i+2}c_{1}c_{2}\cdots c_{2i}c_{2i+1}
\]
as
\[c_{2h}^{-1}\cdots c_{2i+1}^{-1} c_{2h+1}c_{2h}\cdots
c_{2i+2}c_{2i}^{-1}\cdots c_1^{-1}c_{1}c_{2}\cdots c_{2i}c_{2i+1}
\]
\[c_{2h}^{-1}\cdots c_{2i+1}^{-1} c_{2h+1}c_{2h}\cdots
c_{2i+2}c_{2i+1}
\]

Therefore what we had originally is equal to
\[
c_{2i+2}\cdots c_{2h}c_{2h+1}c_{2i}\cdots c_{2}c_{1}c_1\cdots
c_{2h}c_{2h+1}c_{2h}^{-1}\cdots c_{2i+1}^{-1} c_{2h+1}c_{2h}\cdots
c_{2i+2}c_{2i+1}.
\]

Another relation we have to remember is the braid relation
$c_ic_{i+1}c_i=c_{i+1}c_ic_{i+1}$, from which we can obtain
$c_ic_{i+1}c_i^{-1}=c_{i+1}^{-1}c_ic_{i+1}$. Using this multiple
times on the expression
\[c_1\cdots
c_{2h}c_{2h+1}c_{2h}^{-1}\cdots c_{2i+1}^{-1}
\]
along with the commutativity relation mentioned above, we obtain
\[c_{2h+1}^{-1}\cdots c_{2i+2}^{-1}c_1\cdots c_{2h}c_{2h+1}.
\]
Substituting this in what we have for the original expression now
we get
\[ c_{2i+2}\cdots c_{2h}c_{2h+1}c_{2i}\cdots
c_{2}c_{1}c_1\cdots c_{2h}c_{2h+1}c_{2h}^{-1}\cdots c_{2i+1}^{-1}
c_{2h+1}c_{2h}\cdots c_{2i+2}c_{2i+1}
\]
\[
c_{2i+2}\cdots c_{2h}c_{2h+1}c_{2i}\cdots
c_{2}c_{1}c_{2h+1}^{-1}\cdots c_{2i+2}^{-1}c_1\cdots
c_{2h}c_{2h+1} c_{2h+1}c_{2h}\cdots c_{2i+2}c_{2i+1}.
\]
Using the commutativity relation between the terms $c_{2i}\cdots
c_{2}c_{1}$ and $c_{2h+1}^{-1}\cdots c_{2i+2}^{-1}$ we can write
the above expression as
\[
c_{2i+2}\cdots c_{2h}c_{2h+1}c_{2h+1}^{-1}\cdots
c_{2i+2}^{-1}c_{2i}\cdots c_{2}c_{1}c_1\cdots c_{2h}c_{2h+1}
c_{2h+1}c_{2h}\cdots c_{2i+2}c_{2i+1},
\]
which simplifies to
\[c_{2i}\cdots c_{2}c_{1}c_1\cdots c_{2h}c_{2h+1}
c_{2h+1}c_{2h}\cdots c_{2i+2}c_{2i+1}.
\]
If we square this we get
\[c_{2i}\cdots c_{1}\underline{c_1\cdots c_{2h+1}
c_{2h+1}\cdots c_{2i+1}c_{2i}\cdots c_{1}}c_1\cdots c_{2h+1}
c_{2h+1}\cdots c_{2i+1}.
\]
The underlined portion is the well-known expression for $i$. Also
using the fact that $i$ commutes with $c_i,$ the above expression
becomes
\[ic_{2i}\cdots c_{1}c_1\cdots c_{2h+1}
c_{2h+1}\cdots c_{2i+1}.
\]
Now the question reduces to showing
\[ic_{2i}\cdots c_{1}c_1\cdots c_{2h+1}
c_{2h+1}\cdots c_{2i+1}=1.
\]
We will obtain that result by going backwards from the relation
$i^2=1$, by first writing it as
\[ic_1\cdots c_{2h+1}c_{2h+1}\cdots c_{2i+1}c_{2i}\cdots c_{1}=1,
\]
then multiplying by $c_1^{-1}$ on the right,
\[ic_1\cdots c_{2h+1}c_{2h+1}\cdots c_{2i+1}c_{2i}\cdots c_2=c_{1}^{-1},
\]
and then multiplying by $c_1$ on the left
\[ic_1c_1\cdots c_{2h+1}c_{2h+1}\cdots c_{2i+1}c_{2i}\cdots c_2=1
\]
and repeating the same procedure $2i$ times.\\

\begin{figure}[htbp]
     \centering  \leavevmode
     \psfig{file=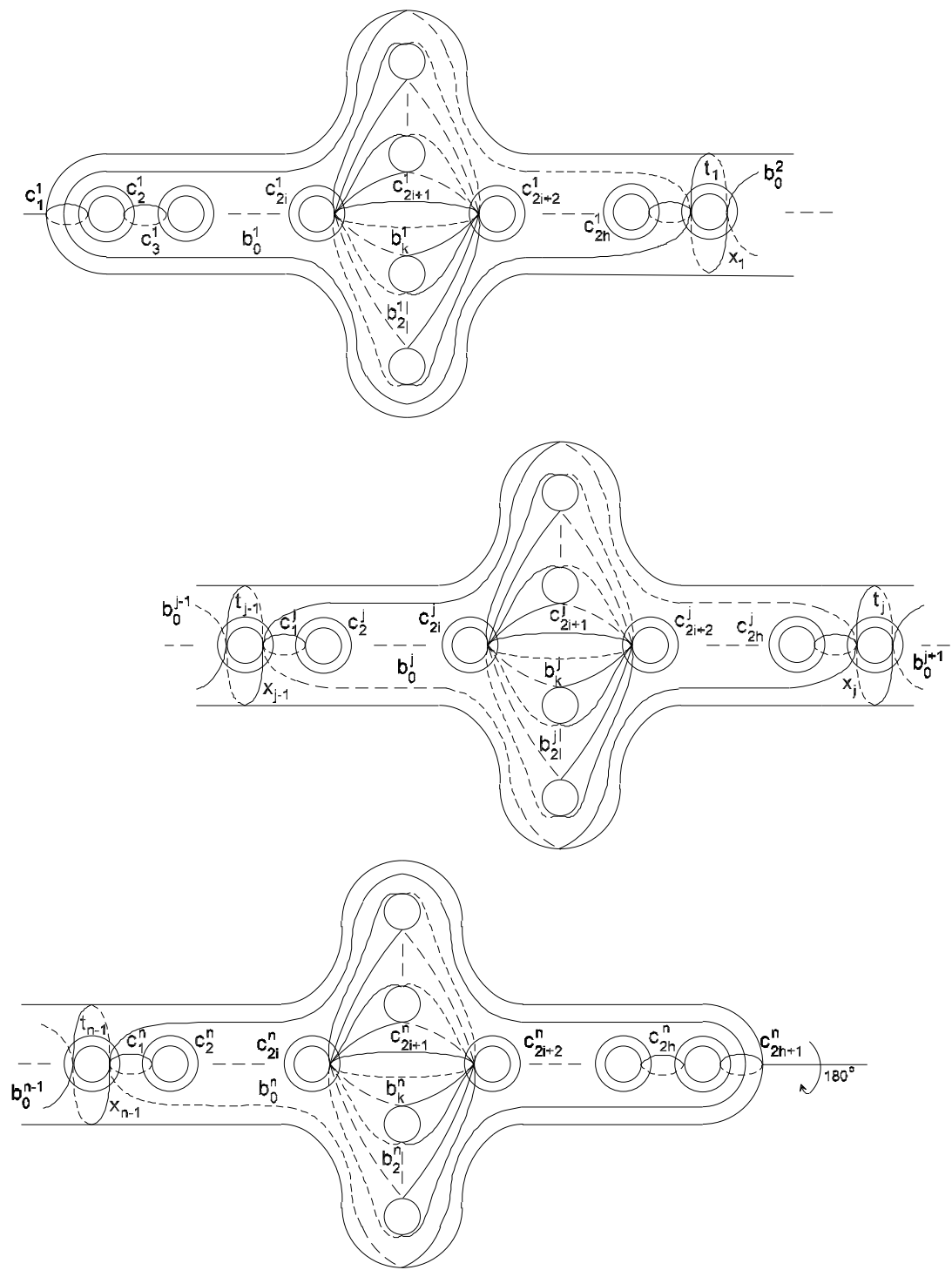,width=14.0in,clip=}
     \caption{An alternate expression for $\theta$}
     \label{altmultiple.fig}
 \end{figure}

\begin{figure}[htbp]
     \centering  \leavevmode
     \psfig{file=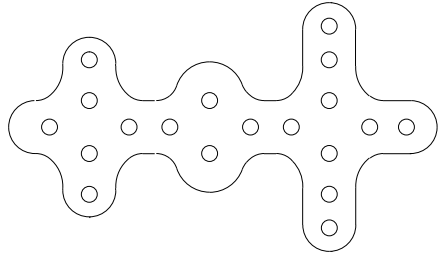,width=2.50in,clip=}
      \caption{A surface with input (1 4 1,1 2 1,1 6 2)}
    \label{inputex1.fig}
 \end{figure}

\begin{figure}[htbp]
     \centering  \leavevmode
     \psfig{file=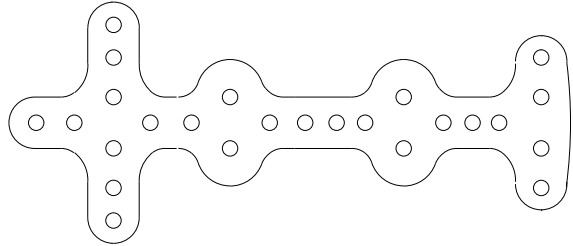,width=3.0in,clip=}
      \caption{A surface with one possible input (2 6 1,1 2 2,2 2 1,2 4 0)}
    \label{inputex2.fig}
 \end{figure}

An alternate expression for $\theta$ using a slightly different
set of cycles is obtained by gluing $n$ copies of bounded surfaces
in Figure \ref{basecycleswithboundary.fig} together along tori
with two boundary components. Figure \ref{altmultiple.fig}
demonstrates the set of cycles that are used in that expression.
The need for this expression emerges from the fact that it is
necessary to have at least two holes between two copies when they
are glued along four-holed spheres, as seen in Figures
\ref{inputex1.fig} and \ref{inputex2.fig}. The alternate
expression allows us to have only one hole between two adjacent
copies and it is very similar to the one given in Theorem
\ref{main.thm}:
\[\theta
=Y^n_i\prod_{j=2}^{n}\left(Y^{j-1}_iX_{j-1}b^j_0Y^j_o\right)b^1_0Y^1_o,
\]

where
\begin{eqnarray*}
   Y^j_o&=&r^j_ol^j_om^j,\\
   Y^j_i&=&r^j_il^j_i,\\
   X_j&=&x_jt_j,
 \end{eqnarray*}
and
 \begin{eqnarray*}
 l^j_i&=&c^j_{2i}\cdots c^j_1,\\ l^j_o&=&c^j_{1}\cdots
c^j_{2i},\\ r^j_i&=&c^j_{2i+2}\cdots c^j_{2h+1},
\\ r^j_o&=&c^j_{2h+1}\cdots c^j_{2i+2}, \\ m^j&=&b^j_1\cdots
b^j_{k}c^j_{2i+1}.
\end{eqnarray*}

\begin{figure}[htbp]
     \centering  \leavevmode
     \psfig{file=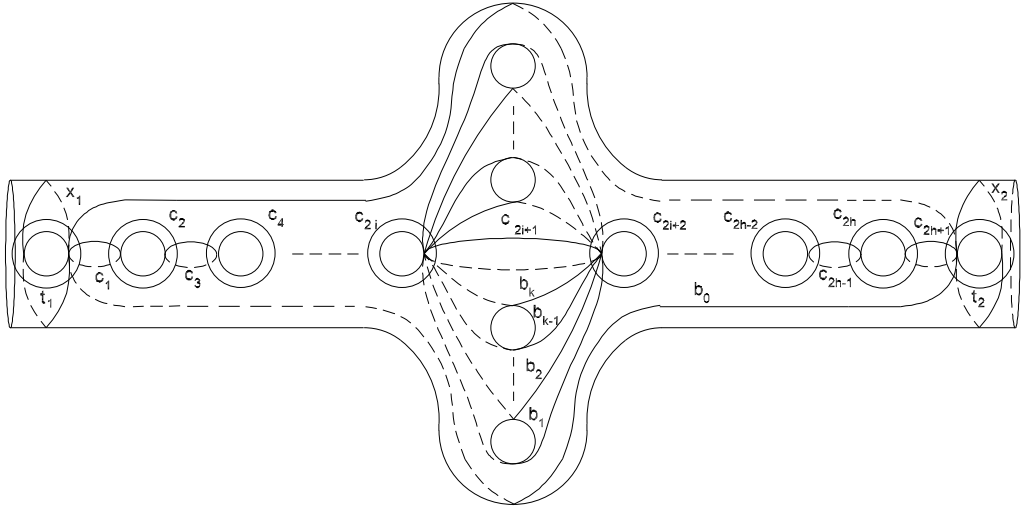,width=5.50in,clip=}
     \caption{Alternate cycles for $\theta$ defined on the bounded surface in Figure \ref{basecycleswithboundary.fig}}
     \label{altsimplecasewithboundary.fig}
 \end{figure}

The proof of this fact also uses induction and it is essentially
based on the simple bounded case that is similar to the one in
Proposition \ref{boundedcase.prop}. After modifying Figure
\ref{basecycleswithboundary.fig} slightly, we obtain Figure
\ref{altsimplecasewithboundary.fig}

and hence the expression
\[
\theta =c_{2i+2}\cdots c_{2h+1}c_{2i}\cdots
c_{1}x_{2}t_{2}x_{1}t_{1}b_{0}c_{2h+1}\cdots c_{2i+2}c_{1}\cdots
c_{2i}b_{1}b_{2}\cdots b_{k-1}b_{k}c_{2i+1}
\]
that replaces the one in Proposition \ref{boundedcase.prop}. We
will not give a detailed proof for this last expression, instead
just provide the mapping of the boundary component $\delta_1$ in
Figure \ref{basecycleswithboundary.fig}. One has to mimic the
steps that is involved in the mapping of the other cycles in the
proof of Proposition \ref{boundedcase.prop} by accommodating the
slight
modifications as needed.\\

\begin{figure}[htbp]
     \centering  \leavevmode
     \psfig{file=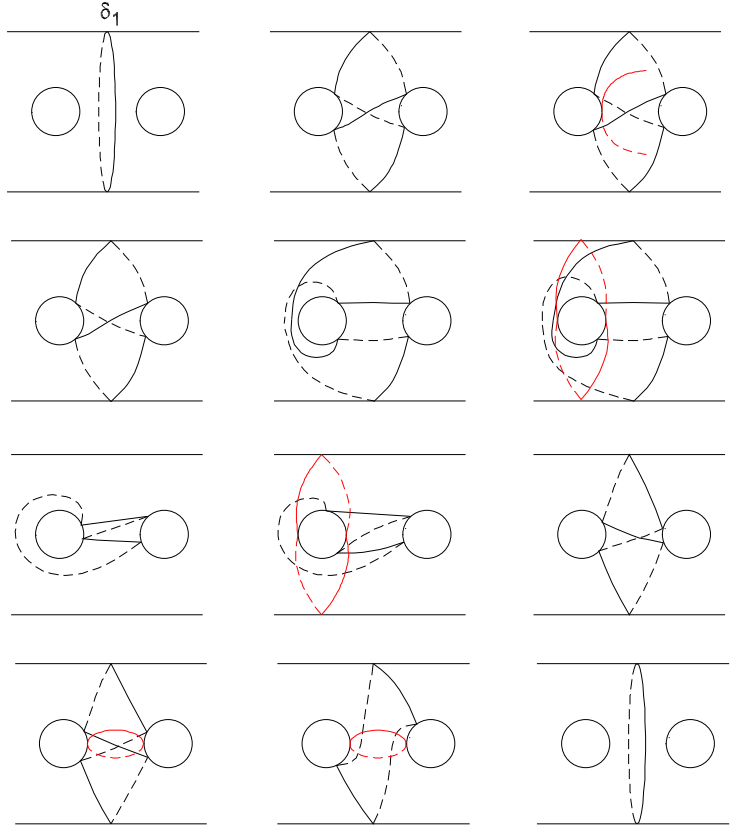,width=8.0in,clip=}
      \caption{The mapping of $\delta_1$}
    \label{altmappingofdelta1.fig}
 \end{figure}

In the first line of Figure \ref{altmappingofdelta1.fig} we see
the effect of $c_1$ on $\delta_1$ first because all the cycles
that come before  $c_1$ miss $\delta_1$. The next cycle, $b_0,$
misses the result from that as seen in the end of the first and
the beginning of the second line. After that, the twist about
$t_1$ takes place, which is not demonstrated in two steps, even
though it intersects the cycle it twists twice. The result from
that has intersection number 2 with $x_1$ and the twist about
$x_1$ is shown in two steps in the end of the second line and all
of the third line. Following twists miss completely the cycle that
is in the end of the third line and the twist about $c_1$ brings
that cycle back to $\delta_1$.\\

A corollary to the expression for $\theta$ in Theorem
\ref{main.thm} and its alternate form would be setting $k=0$ to
obtain some new expressions for the hyperelliptic involution. All
we have to do is redefine $Y^j_o$ as $Y^j_o=r^j_ol^j_oc^j_{2i+1}$
without changing the expression
\[
Y^n_i\prod_{j=2}^{n}\left(Y^{j-1}_iX_{j-1}b^j_0Y^j_o\right)b^1_0Y^1_o,
\]
for $\theta$.
 \clearpage

\section{Applications} In this section we will determine the homeomorphism
type of the genus $g$ Lefschetz fibration
\[
 X\longrightarrow S^{2}
\]
 described by the word $\theta^2 =1$ in the mapping class group
 $M_g$, where $\theta$ is as defined in
Theorem \ref{main.thm}.

Consider the surface in Figure \ref{multiple.fig}. Let $k_j$ be
the $j^{th}$ \emph{vertical genus}, the total genus of the central
part of the $j^{th}$ copy, and let $h_j=l_j+r_j$ be the $j^{th}$
\emph{horizontal genus}, namely the sum of the $j^{th}$ \emph{left
genus} $l_j$ and the $j^{th}$ \emph{right genus} $r_j$. Let
$k=\sum k_j$ be  \emph{the vertical genus} and $h=\sum h_j$ be
\emph{the horizontal genus}. If we denote the total genus by $g$
then $g=h+k$.

To find the total number of cycles contained in $\theta$ let's
recall that

\[\theta
=Y^n_i\prod_{j=2}^{n}\left(Y^{j-1}_iX_{j-1}b^j_0Y^j_o\right)b^1_0Y^1_o,
\]
 where
 \begin{eqnarray*}
   Y^j_o&=&r^j_ol^j_om^j,\\
   Y^j_i&=&r^j_il^j_i,\\
   X_j&=&x_jt_j,
 \end{eqnarray*}

and

 \begin{eqnarray*}
 l^j_i&=&c^j_{2i}\cdots c^j_2,\\ l^j_o&=&c^j_{2}\cdots
c^j_{2i},\\ r^j_i&=&c^j_{2i+2}\cdots c^j_{2h},
\\ r^j_o&=&c^j_{2h}\cdots c^j_{2i+2}, \\ m^j&=&b^j_1\cdots
b^j_{k}c^j_{2i+1}.
\end{eqnarray*}

The first and the last copies would differ in the first two and
the following two lines of definitions above, respectively:

\begin{eqnarray*}
 l^1_i&=&c^1_{2i}\cdots c^1_1,\\ l^1_o&=&c^1_{1}\cdots
c^1_{2i},\\ r^n_i&=&c^n_{2i+2}\cdots c^n_{2h+1},
\\ r^n_o&=&c^n_{2h+1}\cdots c^n_{2i+2}.
\end{eqnarray*}

$m^j$ consists of $k_j+1$ cycles. Both $l^j_i$ and $l^j_o$ consist
of $2i-2+1=2i-1$ cycles for $j\neq1$ and $l^1_i$ and $l^1_o$
consist of $2i-1+1=2i$ cycles. Likewise, both $r^j_i$ and $r^j_o$
consist of $2h_j-(2i+2)+1=2h_j-2i-1$ cycles for $j\neq n$ and
$r^n_i$ and $r^n_o$ consist of $2h_n+1-(2i+2)+1=2h_n-2i$ cycles.
Therefore $Y^j_o$ consists of
$y^j_o=2h_j-2i-1+2i-1+k_j+1=k_j+2h_j-1$ cycles for $j\neq1,n$,
$Y^1_o$ consists of $y^1_o=2h_1-2i-1+2i+k_1+1=k_1+2h_1$ cycles,
and $Y^n_o$ consists of $y^n_o=2h_n-2i+2i-1+k_n+1=k_n+2h_n$
cycles. $Y^j_i$ has $y^j_i=2h_j-2i-1+2i-1=2h_j-2$ cycles for
$j\neq1,n$, $Y^1_i$ has $y^1_i=2h_1-2i-1+2i=2h_1-1$ cycles, and
$Y^n_i$ has $y^n_i=2h_n-2i+2i-1=2h_n-1$ cycles. Clearly $X_j$
consists of 2 cycles.

In the above computations, too, we ignored the dependence of $i$
on $j$ and did not write $2i_j$ instead of $2i$ in order not to
make the computations more complicated because they cancel out
anyway.

Now, using the lengths of each group of twists computed above we
determine that

\[Y^n_i\prod_{j=2}^{n}\left(Y^{j-1}_iX_{j-1}b^j_0Y^j_o\right)b^1_0Y^1_o,
\]

consists of

\[y^n_i+\sum_{j=2}^{n}\left(y^{j-1}_i+2+1+y^j_o\right)+1+y^1_o
\]

\[=y^n_i+\sum_{j=2}^{n}
y^{j-1}_i+3(n-2+1)+\sum_{j=2}^{n}y^j_o+1+y^1_o
\]

cycles. Rearranging the indices and simplifying we get

\[y^n_i+\sum_{j=1}^{n-1}y^{j}_i+3(n-1)+\sum_{j=2}^{n}y^j_o+1+y^1_o
\]

Releasing the first term of the first sum and the last term of the
second sum gives

\[y^n_i+y^1_i+\sum_{j=2}^{n-1}y^{j}_i+3(n-1)+y^n_o+\sum_{j=2}^{n-1}y^j_o+1+y^1_o,
\]

which is equal to

\[y^1_i+y^1_o+y^n_i+y^n_o+\sum_{j=2}^{n-1}\left(y^{j}_i+y^j_o\right)+3(n-1)+1.
\]

Now substituting the value of each term and simplifying we obtain

\[2h_1-1+k_1+2h_1+2h_n-1+k_n+2h_n+\sum_{j=2}^{n-1}\left(2h_j-2+k_j+2h_j-1\right)+3(n-1)+1
\]

\[=4h_1+k_1+4h_n+k_n+\sum_{j=2}^{n-1}\left(4h_j+k_j\right)-3(n-1-2+1)+3(n-1)-1
\]

\[=\sum_{j=1}^{n}\left(4h_j+k_j\right)-3(n-2)+3(n-1)-1
\]

\[=4\sum_{j=1}^{n}h_j+\sum_{j=1}^{n}k_j+6-3n+3n-3-1
\]

\[=4h+k+2
\]

Therefore the word $\theta^2 =1$ consists of $2(4h+k+2)=8h+2k+4$
twists.

Since all the twists are about non-separating cycles, the
Lefschetz fibration defined by the word $\theta^2 =1$  has
$8h+2k+4$ irreducible fibers. This allows us to compute the Euler
characteristic of the $4-$ manifold $X$ using the formula

\[ \chi(X)=2(2-2g)+\mbox{number of singular fibers}
\]

for Lefschetz fibrations, which is

 \[2(2-2g)+8h+2k+4=4-4g+8h+2k+4=4-4(h+k)+8h+2k+4
 \]

 \[=8+4h-2k
 \]

 in our case.

 The other homeomorphism invariant that we will compute is the
 signature $\sigma(X)$ of the $4-$ manifold $X$.

 Using the algorithm described in \cite{Oz} we wrote a Matlab
program that computes the signature of the Lefschetz fibration
described by the word $\theta^2 =1$.

The input for the program is the left, right, and the vertical
genus of each copy that is glued together to form the surface
$\Sigma$ on which $\theta$ is defined. The following are two
examples that demonstrate how the shape of the surface is coded
into a sequence of numbers, which are used as the inputs for the
program.

\clearpage (0 2 1,1 2 0):

$0+0+0+0-1+0+0+0+0-1-1-1-1-1-1-1+0-1+0-1-1+0-1+0+0+0+0+0=-12$

\bigskip

(0 4 1,1 2 0):

$0+0+0+0+0+0-1+0+0+0+0-1-1-1-1-1-1-1+0+0+0-1+0-1-1+0-1+0+0+0+0+0=-12$

\bigskip

(0 2 1,1 4 0):

$0+0+0+0-1+0+0+0+0+0+0-1-1-1-1-1-1-1+0-1+0-1-1+0+0+0-1+0+0+0+0+0=-12$

\bigskip

(0 4 1,1 4 0):

$0+0+0+0+0+0-1+0+0+0+0+0+0-1-1-1-1-1-1-1+0+0+0-1+0-1-1+0+0+0-1+0+0+0+0+0=-12$

\bigskip

(1 2 1,1 4 0):

$0+0+0+0+0+0-1+0+0+0+0+0+0-1-1-1-1-1-1-1-1-1+0-1-1-1+0-1-1+0+0+0-1+
0+0+0+0+0+0+0=-16$

\bigskip

(0 2 1,1 4 1):

$0+0+0+0-1+0+0+0+0+0+0+0+0-1-1-1-1-1-1-1-1-1+0-1+0-1-1+0+0+0-1-1-1+
0+0+0+0+0+0+0=-16$

\bigskip

(0 2 2,1 4 0):

$0+0+0+0+0+0-1+0+0+0+0+0+0-1-1-1-1-1-1-1-1-1+0-1-1-1+0-1-1+0+0+0-1+
0+0+0+0+0+0+0=-16$

\bigskip

(0 2 1,1 2 1,1 2 0):

$0+0+0+0-1+0+0+0+0+0+0-1-1-1+0+0+0+0-1-1-1-1-1-1-1-1+0-1+0-1-1+0+
-1-1+0+0+0+0-1-1+0-1+0+0+0+0+0+0=-20$

\bigskip

(0 2 1,1 2 1,1 2 2,1 2 1):

$0+0+0+0-1+0+0+0+0+0+0-1-1-1+0+0+0+0+0+0+0+0-1-1-1-1+0+0+0+0+0+0-1
-1-1-1-1-1-1-1-1-1-1-1+0-1+0-1-1+0-1-1+0+0+0+0-1-1+0-1-1-1-1+0+0+0+0+0-1-1
+0-1-1-1+0+0+0+0+0+0+0+0+0+0=-36$

\bigskip

(2 2 1,1 2 2,1 4 1):

$0+0+0+0+0+0+0+0-1+0+0+0+0+0+0+0+0-1-1-1-1-1-1-1+0+0+0+0+0+0+0+0-1-1-1-1-1-1-1
-1-1-1-1-1+0-1-1-1-1-1+0-1-1+0-1-1-1-1+0+0+0+0+0+0+0+0-1-1+0+0+0-1-1-1+0+0+0+0+0+0+0+0+0+0=-36$

\bigskip

(3 4 2,1 4 2):

$0+0+0+0+0+0+0+0+0+0+0+0+ 0+0-1+0+0+0+0+0+0+0+0+0+0-1-1-1-1-1-1
-1-1-1-1-1-1-1-1-1-1-1-1-1+0+0+0-1-1-1-1-1-1-1-1-1+0-1-1+0+0+0-1-1-1-1-1
+0+0+0+0+0+0+0+0+0+0+0+0+0+0+0+0+0=-36$

\bigskip

(1 4 1,1 2 1,1 6 2): -32

\clearpage

The above computations, along with many others that we do not
include here, point out to the fact that the signature depends
only on $h$, i.e., it is independent of $k$. A quick check
suggests that $\sigma(X)=-4(h+1)$ for the above computations. We
conjecture that this is true in general, namely the signature of
the Lefschetz fibration given by the word $\theta^2 =1$, where
$\theta$ is as defined in Theorem \ref{main.thm}, is $-4(h+1)$.

For $\chi(X)=8+4h-2k$ and $\sigma(X)=-4(h+1)$, we obtain

\begin{eqnarray*}
c_1^2(X)&=&3\sigma(X)+2\chi(X) \\
&=&3(-4h-4)+2(8+4h-2k) \\
&=& -4h-4k+4 \\
&=& -4(g-1)
\end{eqnarray*}

and

\begin{eqnarray*}
\chi_h(X)&=&\frac{1}{4}(\sigma(X)+\chi(X)) \\
&=&\frac{1}{4}(8+4h-2k-4h-4) \\
&=& \frac{1}{4}(4-2k)=1-k/2
\end{eqnarray*}

Recall that $k$ is even.  $\chi_h(X)$ makes sense here because $X$
has  almost complex structure.

\bibliographystyle{amsplain}

\end{document}